\documentclass{amsart}

\title{Desperately seeking mathematical truth}
\author{Melvyn B. Nathanson}
\address{Department of Mathematics,
Lehman College (CUNY),
Bronx, New York 10468, and 
CUNY Graduate Center, New York, NY 10016}
\email{melvyn.nathanson@lehman.cuny.edu}

\begin{document}
\maketitle
 
We mathematicians have a naive belief in truth.  We prove theorems.  Theorems are deductions from axioms.   Each line in a proof is a simple consequence of the previous lines of the proof, or of previously proved theorems.  Our conclusions are true, unconditionally and eternally.  The Babylonians' quadratic formula and the Greeks' proof of the irrationality of $\sqrt{2}$ are true even in the Large Magellanic Cloud.  

How do we know that a proof is correct?  By checking it, line by line.  A computer might even be programmed to check it.  To discover a proof, to have the imagination and genius to conceive a chain of reasoning that leads from trivial axioms to extraordinarily beautiful conclusions--this is a rare and wonderful talent.  This is mathematics.  But to check a proof--any fool can do this.  

Still, there is a nagging worry about this belief in mathematical certitude.  In 2000 the Clay Mathematics Institute announced million dollar prizes for the solution of seven ``Millennium problems.''  Solve one of the problems and receive a million dollars.  According to CMI's rules, two years after the appearance of the solution in a ``refereed mathematics publication of worldwide repute'' and after ``general acceptance in the mathematics community,'' the prize would be awarded.  

But why the delay?  Surely, any competent person can check a proof.  It's either right or wrong.  Why wait two years?

The reason is that many great and important theorems don't actually have proofs.  They have sketches of proofs, outlines of arguments, hints and intuitions that were obvious to the author (at least, at the time of writing) and that, hopefully, are understood and believed by some part of the mathematical community. 

But the ÒcommunityÓ itself is tiny.  In most fields of mathematics there are few experts.  Indeed, there are very few active research mathematicians in the world, and many important problems, so the ratio of the number of mathematicians to the number of problems is small.  In every field, there are  ``bosses'' who proclaim the correctness or incorrectness of a new result, and its importance or unimportance. Sometimes they disagree, like gang leaders fighting over turf.  In any case, there is a web of semi-proved theorems throughout mathematics.  Our knowledge of the truth of a theorem depends on the correctness of its proof and on the correctness of all of the theorems used in its proof.  It is a shaky foundation.  

Even Euclid got things wrong, in the sense that there are statements in the Elements (e.g. Book I, Proposition 1) that do not follow logically from the axioms.  It took 150 years after Leibnitz and Newton until the foundations of differential and integral calculus were formulated correctly, and we could backfill proofs of much 18th and 19th century mathematical analysis.

Similar problems plague us today.  Consider the two most highly publicized recent successes of mathematics.  It took several years to confirm the correctness of Wiles' proof of ``Fermat's Last Theorem.''  A mistake was found in the original paper, and there still remained questions about the truth of other results used in the proof.  There were also arguments about the completeness of Perelman's proof of the Poincar{\' e} conjecture, the first of the Millennium problems to be solved.  How many mathematicians have checked both Wiles' and Perelman's proofs?  

I certainly don't claim that there are gaps in Wiles' or Perelman's work.  I don't know.  We (the mathematical community) believe that the proofs are correct because a political consensus has developed in support of their correctness.  

The classification of the finite simple groups provides another example.  There is still uncertainty about whether there is a complete proof of the classification, and, if there is a proof, exactly when the theorem was proved.   After the classification was supposedly finished, Danny Gorenstein wanted to write an exposition of the proof.  He would reread the papers and find mistakes.  He could always patch things up, but until he did, the ÒtheoremsÓ in the original papers were, according to the rules of our profession, not theorems but unproven assertions.  We mathematicians like to talk about the ``reliability'' of our literature, but it is, in fact, unreliable.

Part of the problem is refereeing.  Many (I think most) papers in most refereed journals are not refereed.  There is a presumptive referee who looks at the paper, reads the introduction and the statements of the results, glances at the proofs, and, if everything seems OK, recommends publication.  Some referees do check proofs line-by-line, but many do not.  When I read a journal article, I often find mistakes.  Whether I can fix them is irrelevant.  The literature is unreliable. 

How do we recognize mathematical truth?  If a theorem has a short complete proof, we can check it.  But if the proof is deep, difficult, and already fills100 journal pages, if no one has the time and energy to fill in the details, if a  ``complete'' proof would be 100,000 pages long, then we rely on the judgments of the bosses in the field.  In mathematics, a theorem is true, or it's not a theorem.  But even in mathematics, truth can be political.  

\end{document}